\documentclass[12pt,leqno]{article}
\usepackage{amsmath,amssymb,latexsym}
\usepackage{amscd}
\usepackage{mathrsfs}
\usepackage{srcltx}
\usepackage{theorem}
\usepackage{amsfonts}

\newcommand{\R}{\mathbb{R}}        
\newcommand{\N}{\mathbb{N}}        
\newcommand{\Rn}{\mathbb{R}^d}     

\newcommand{\CN}{C^\infty(\Rn)}

\newcommand{\OH}{\cO_H'(\Rn)}
\newcommand{\Ohnz}{\cO_H(\NZ)}
\newcommand{\OHNZ}{\cO_H'(\NZ)}
\newcommand{\OHq}{\cO_H'(Q)}

\newcommand{\OY}{\cO_Y'(\Rn)}
\newcommand{\Oy}{\cO_Y(\Rn)}

\newcommand{\Di}{\cD(\Rn)}
\newcommand{\di}{\cD(\R)}
\newcommand{\Dip}{\cD\,'(\Rn)}
\newcommand{\Si}{\cS(\Rn)}
\newcommand{\Sip}{\cS\,'(\Rn)}
\newcommand{\Siq}{\cS(Q)}
\newcommand{\Sipq}{\cS\,'(Q)}
\newcommand{\si}{\cS(\R)}
\newcommand{\sip}{\cS\,'(\R)}

\newcommand{\NZ}{\Rn_*}

\newcommand{\Proof}{\textbf{Proof:} \ }
\newcommand{\qed}{\hspace*{\fill} $\Box $}
\newcommand{\To}{\longrightarrow}
\newcommand{\cA}{\mathscr{A}}
\newcommand{\cB}{\mathscr{B}}
\newcommand{\cF}{\mathcal{F}}
\newcommand{\cD}{\mathscr{D}}

\newcommand{\cS}{\mathscr{S}}
\newcommand{\cO}{\mathscr{O}}  

\newcommand{\be}{{\bf 1}}

\newcommand{\F}{Fr\'echet }
\renewcommand{\span}{\mathrm{span}}

\newcommand{\om}{\omega}

\newtheorem{proposition}{Proposition}[section]

\newtheorem{lemma}[proposition]{Lemma}
\newtheorem{corollary}[proposition]{Corollary}
\newtheorem{theorem}[proposition]{Theorem}

\newtheorem{definition}{Definition}

\newcommand{\Exp}{{\rm Exp}}
\newcommand{\Log}{{\rm Log}}

\newcommand{\supp}{\mathrm{supp}\, }

\newcommand{\vp}{\varphi}

\newcommand{\tL}{\widetilde{L}}

\title{{\sc Hadamard type operators on temperate distributions
}}
\author{Dietmar Vogt}
\date{}

\begin{document}

\maketitle

\footnotetext{\hskip -.8em   {\em 2010 Mathematics Subject
Classification.}
    {Primary: 46F10. Secondary: 47B38,44A35,46E10.}
    \hfil\break\indent \begin{minipage}[t]{14cm}{\em Key words and phrases:} Hadamard operators, temperate distributions, convolution operators, $C^\infty$-functions of exponential decay. \end{minipage}.
     \hfil\break\indent
{}}

\begin{abstract} We study Hadamard operators  on $\Sip$ and give a complete characterization. They have the form $L(S)=S\star T$ where $\star$ means the multiplicative convolution and $T\in\OHNZ$, the space of distributions which are $\theta$-rapidly decreasing in infinity and at the coordinate hyperplanes. To show this we study and characterize convolution operators on the space $Y(\Rn)$ of exponentially decreasing $C^\infty$-functions on $\Rn$. We use this and the exponential transformation to characterize the Hadamard operators on $\Sipq$, $Q$ the positive quadrant, and this result we use as a building block for our main result.

\end{abstract}



In the present note we study Hadamard operators on $\Sip$, that is, continuous linear operators on $\Sip$ which admit all monomials as eigenvectors and we give a complete characterization. Operators of Hadamard type have attracted some attention in recent times. Such operators on $\CN$ have been studied and characterized in \cite{Vhad, Vconv}, on $\cA(\R)$ in \cite{DLre, DLalg, DLhad} and on $\cA(\Rn)$ in \cite{DLV}. There you find also references to the long history of such problems.  Their surjectivity on $\CN$ has been characterized in \cite{DLhadsurj}. Since it can be shown that Hadamard operators commute with dilations our problem is, by duality, closely related to the study of continuous linear operators in $\Si$ which commute with dilations. In a first step we study such operators on $\Siq$,  $Q=]0,+\infty[^d$. By means of the exponential transformation this can be transferred to the study of convolution operators on the space $Y(\Rn)$ of $C^\infty$-functions on $\Rn$ with exponential decay.

In a first part of the paper we study such operators and give a complete characterization in terms of the class $\OY$ of exponentially decreasing distributions, which is similar to the class $\cO_C'$ of L. Schwartz of rapidly decreasing distributions, which are the convolution-multipliers in $\Si$. We study the class $\OY$ and these results are of independent interest.

By means of the exponential transformation we obtain a description of the operators on $\Sipq$ which commute with dilations in $Q$.
They have the form $\vp\mapsto T_x\vp(xy)$ where $T$ is a distribution in $\OHq$. These are the exponential transforms of $\OY$, we call them $\theta$-rapidly decreasing distributions on $Q$. The class $\OH$ first appeared in \cite{VDprime1} where the Hadamard operators in $\Dip$ were described.  For a more detailed study of this class and examples see  \cite[\S3]{VDprime1}.

From there we obtain our main result: The Hadamard operators on $\Sip$ have the form $S\mapsto S\star T$ where $T\in\OHNZ$ the class of distributions on $\Rn$ which are $\theta$-rapidly decreasing in infinity and at the coordinate hyperplanes. It is a subclass of $\OH$, known from \cite{VDprime1}.

We use standard notation of Functional Analysis, in  particular, of distribution theory. For unexplained notation we refer to \cite{DK}, \cite{LSI}, \cite{LS}, \cite{MV}.

\section{Preliminaries} \label{s0}
We use the following notation $\partial_j=\partial/\partial x_j$, $\theta_j = x_j\,\partial_j$. For a multiindex $\alpha\in\N_0^d$ we set $\partial^\alpha = \partial_1^{\alpha_1}..\partial_d^{\alpha_d}$, likewise for $\theta^\alpha$. $\be$ denotes the vector $(1,\dots,1)$.
For vectors $x,y\in\Rn$ we will use the definition $xy=(x_1 y_1,...x_d y_d)$. This will hold except for obvious cases like in the formula for the Fourier transform.

For a polynomial $P(z)=\sum_\alpha c_\alpha z^\alpha$ we consider the  \it Euler operator \rm $P(\theta)=\sum_\alpha c_\alpha \theta^\alpha$ and also the operator $P(\partial)$, defined likewise. The dual operator of $P(\theta)$ is $P(\theta^*)$ where $\theta^*=-\theta-1$, hence also an Euler operator.

For $a\in\NZ$ the \em dilation operator \rm $D_a$ is defined by $(D_a T)\vp= |a_1\cdots a_d|^{-1} T_\xi\vp(\xi/a)$. For the distribution $x^\alpha\in\Sip$ this yields $D_a x^\alpha= (ax)^\alpha$. For $e\in\{-1,+1\}^d$ this definition simplifies to
$(D_e T)\vp =T_\xi\vp(e\xi)$. These operators are called \em reflections \rm.

For basic properties of Hadamard operators see \cite{VDprime1}. They are a closed commutative sub-algebra of $L(\Sip)$. Euler operators and dilations are of Hadamard type, Therefore they commute with all Hadamard operators. On the other hand we have:

\begin{lemma}\label{l9} If $L\in L(\Sip)$ commutes with $\theta_j$ for all $j$ and with all reflections then it is a Hadamard operator.
\end{lemma}

\Proof We set $T=L(x^\alpha)$ and have to show that $T\in\span\{x^\alpha\}$. Since $L$ commutes with $\theta_j$ we obtain $\theta_j T=\alpha_j T$. By use of the exponential transformation we obtain for $Q$ and likewise for all quadrants $Q_e=eQ$ that $T=c_e x^\alpha$ on $Q_e$, with constants $c_e$. Since $L$ commutes with reflections all $c_e$ must be equal and we have $T=c x^\alpha$ on $\NZ$.

We set $S=T-c x^\alpha$. Then $\supp S\subset Z_0=\{\xi\,:\,\xi_1\cdots\xi_d=0\}$ and $\theta_j S= \alpha_j S$. Since $S$ is of finite order there is $\beta\in\N^d$ such that $x^\beta S=0$. We have
$$\partial_j(x^\beta S)=\beta_j x^{\beta'} S+ x^{\beta'}\theta_j S=(\beta_j+\alpha_j)\, x^{\beta'} S$$
where $\beta'=(\beta_1,\dots,\beta_j-1,\dots,\beta_d)$. Repeating this we obtain:
$$0=\partial^\beta (x^\beta S)= b\, S$$
with $b\neq 0$. Therefore $S=0$, that is, $L(x^\alpha)=c x^\alpha$. \qed

We set for $x\in\Rn$
$$\Exp(x)=(\exp(x_1),..,\exp(x_d)).$$
$\Exp$ is a diffeomorphism from $\Rn$ onto $Q:=(0,+\infty)^d$. Therefore
$$C_\Exp:f\To f\circ\Exp$$
is a linear topological isomorphism from $C^\infty(Q)$ onto $C^\infty(\Rn)$. For $f\in C^\infty(Q)$ we have
$P(\partial)(f\circ\Exp)=(P(\theta)f)\circ\Exp$ that is $P(\partial)\circ C_\Exp =C_\Exp\circ P(\theta)$. In this way the study of Hadamard operators on $Q$ can be reduced to the study of operators on $\Rn$. This has be done in $\cite{Vconv}$ for $C^\infty(Q)$. We apply the same argument to the space $\cS(Q)$ where $\cS(Q)=\{f\in\cS(\Rn)\,:\,\supp f\subset \overline{Q}\}$.

As usual $\cS(\Rn)$ denotes the Schwartz space of rapidly decreasing $C^\infty$-functions on $\Rn$, its dual $\Sip$   the space of temperate distributions. We consider is subspace $\cS(Q)$ and its dual $\cS'(Q)$.

We recall the following definitions of \cite[Chap. VI, \S 8]{LS}: $\cB'$ denotes the dual of the space of $C^\infty$-space which are bounded including all derivatives
and
$\cD_{L_1}'$ the dual of the space of $C^\infty$-space such that all derivatives are in $L_1(\Rn)$.

\section{Convolution operators on $C^\infty$-functions with exponential decay} \label{s1}

We start with studying convolution operators on the space of $C^\infty$-functions with exponential decay on $\Rn$ and its dual. We will transfer our results by the exponential diffeomorphism to results on Hadamard operators on $\Sipq$  and use this as building blocks to study Hadamard operators on $\Sip$. We set
\begin{eqnarray*} Y(\Rn)&:=&\{f\in C^\infty(\Rn)\,:\,\sup_x|f^{(\alpha)}(x)|\,e^{k|x|}<\infty \text{ for all  }\alpha\text{ and }k\in\N\}\\
&=& \{f\in C^\infty(\Rn)\,:\,\sup_x|f^{(\alpha)}(x)|\,e^{x\eta}<\infty \text{ for all  }\alpha\text{ and }\eta\in\Rn\}
\end{eqnarray*}
with its natural topology.

Then $Y(\Rn)$ is a \F space, closed under convolution and $P(\partial)$  is a continuous linear operator in $Y(\Rn)$ for every polynomial $P$. $\Di\subset Y(\Rn)$ as a dense subspace, hence $Y(\Rn)'\subset \Dip$. We obtain (see \cite[Lemma 2.1]{VEs}):

\begin{lemma} \label{l1} $C_\Exp(\cS(Q))= Y(\Rn)$.
\end{lemma}

We set
$\om(x)=\sum_{\eta\in\{-1,+1\}^d} e^{\eta x}$. We have $\om\in\CN$ and $e^{|x|}\le \om(x)\le 2^d e^{|x|}$.

In analogy to \cite[Chap. VII, \S 5, p. 100]{LS} we define
\begin{definition}\label{d1}
$T\in\OY$ if $\om(kx)T\in \cB'$ for every k.
\end{definition}
It is obvious that we might equivalently write $\om(kx)T\in \cD_{L_1}'$ for every k.

For the following theorem compare \cite[Chap. VII, \S 5, Th\'eor\`eme IX]{LS}.
\begin{theorem}\label{t1} For $T\in\Dip$ the following are equivalent:
\begin{enumerate}
\item $T\in\OY$.
\item For any $k$ there are finitely many functions $t_\beta$ such that $e^{k|x|}t_\beta\in L_\infty(\Rn)$ and such that $T=\sum_\beta \partial^\beta t_\beta$.
\item $T\in Y(\Rn)'$ and $T_x \vp(x+y)\in Y(\Rn)$ for all $\vp\in Y(\Rn)$.
\item $f(y)=T_x \vp(x+y)$ is a exponentially decreasing continuous function (that is $\sup |f(y)| e^{k|y|}<\infty$ for all $k$) for all $\vp\in\Di$.
\item $(\om(kx)T)*\vp$ is a continuous bounded function for every $k$ and $\vp\in\Di$.
\end{enumerate}
\end{theorem}

\Proof (1) $\Rightarrow$ (2) If $\om(kx)T\in\cD_{L_1}'$ then, by a standard conclusion, there are finitely many functions $\tau_\beta\in L_\infty(\Rn)$ such that $\om(kx)T=\sum_\beta \partial^\beta \tau_\beta$. This yields
\begin{eqnarray*}
T\vp &=& (\om(kx) T)\left(\frac{1}{\om(kx)}\,\vp\right)=\sum_\beta \partial^\beta\tau_\beta \left(\frac{1}{\om(kx)}\,\vp\right)\\
&=& \sum_\beta(-1)^{|\beta|}\int \tau_\beta(x) \partial^\beta \left(\frac{1}{\om(kx)}\,\vp(x)\right) dx\\
&=& \sum_\beta(-1)^{|\beta|}\int \tau_\beta(x) \sum_{\alpha\le\beta} c_{\alpha,\beta}\left( \partial^{\beta-\alpha}\frac{1}{\om(kx)}\right)\vp^{(\alpha)}(x) dx\\
&=& \sum_\alpha \left(\sum_{\beta\ge\alpha} (-1)^{|\beta|} \tau_\beta(x) c_{\alpha,\beta}\partial^{\beta-\alpha}\frac{1}{\om(kx)}\right) \vp^{(\alpha)}(x) dx\\
&=& \big\langle \sum_\alpha \partial^\alpha t_\alpha, \vp\big\rangle
\end{eqnarray*}

where $\om(kx)t_\alpha\in L_\infty(\Rn)$ for all the, finitely many, $\alpha$.

(2) $\Rightarrow$ (1) is straightforward, because we may assume that $T=\partial^\beta t_\beta$.

(2) $\Rightarrow$ (3) The first part is clear from (2). Assume $T=\partial^\beta t_\beta$, $e^{(k+1)|x|}|t_\beta(x)|\in L_\infty(\Rn)$. Then
$$T_x\vp(x+y)=(-1)^{|\beta|}\int t_\beta(x) \vp^{(\beta)}(x+y) dx\in \CN$$ and we have
$$e^{k|y|} |\partial_y^\alpha T_x\vp(x+y)|\le \int e^{k|x|} |t_\beta(x)| e^{k|x+y|}|\vp^{(\alpha+\beta)}(x+y)| dx\le \infty.$$
If $T\in\OY$ this holds for all summands in the representation of $T$ with given $k$ and since we have for all $k$ such a representation the claim is proved.

(3) $\Rightarrow$ (4) is obvious.

(4) $\Rightarrow$ (5) For $\vp\in\Di$ and $\eta\in \Rn$ we obtain
$$((e^{\eta x}T)*\vp) (y)= T_x (e^{\eta x} \vp(y-x)) = e^{\eta y} T e^{-\eta(y-x)}\vp(y-x)=e^{\eta y} (T*(e^{-\eta x}\vp))(y).$$
Since $e^{-\eta x}\vp \in\Di$ the right hand side is bounded, by (4). Adding over all $\eta\in \{+k,-k\}^d$ we obtain the result.

(5) $\Rightarrow$ (2) This follows from Lemma \ref{l5}. \qed

The following Lemma is essentially an adaptation of \cite[Chap. VI, \S 8, Th\'eor\`eme XXV]{LS}.

\begin{lemma}\label{l5} Let $\om$ be measurable, $\om(x)> 0$ for all $x\in\Rn$. Let $S\in \Dip$ be a distribution such that
$\sup_x \om(x) |S_y\vp(x-y)|<\infty$ for all $\vp\in \Di$
then there are finitely many measurable functions $\tau_\beta$ with $\sup_x \om(x) \tau_\beta(x)<\infty$ such that $S=\sum_\beta \tau_\beta^{(\beta)}$.
\end{lemma}

\Proof We consider the map $\Psi:\Di\to L_\infty(\Rn)$ given by
$$\Psi(\vp)=\om(x) S_y\vp(x-y).$$
Because of the Closed Graph Theorem $\Psi$ is continuous. Let $B$ denote the unit ball in $\Rn$.
Then there is $m\in\N$ such that $\Psi$ restricted to $\cD(B)$ extends to a continuous map $\cD^{m}(B)\to L_\infty(\Rn)$, where $\cD^{m}(B)$ denotes the Banach space of $m$-times continuously differentiable functions with support in $B$. We choose $\gamma\in\cD(B)$, $\gamma(x)=1$ in a neighborhood of $0$ and set $g=\gamma E\in \cD^m(B)$ where $E$ is an elementary solution of $\Delta^k$, $k$ large enough. Then
$\Psi(g)\in L_\infty(\Rn)$ that means $\tau:=S*g$ is a measurable function with $\om(x) |\tau(x)|\le C$ for suitable $C$ and we obtain $\Delta^k \tau = S+S*\psi $ where $\psi \in \cD(B)$. We have $\om(x) (S*\psi)=\Psi(\psi)\in L_\infty(\Rn)$. Therefore the equality $S=\Delta^k \tau - S*\psi$ shows the result.
\qed

\medskip

We have to fix our notation on the convolution of distributions. For distributions $T, S$ and a function $\psi$ we define $(S*T)\psi :=S_y(T_x \psi(x+y))$ whenever this makes sense.

\begin{lemma}\label{l1} If $T\in\OY$ and $\vp\in Y(\Rn)$ then both $T*\vp\in Y(\Rn)'$ and $\vp*T\in Y(\Rn)'$ are defined and equal and we have $T*\vp=\vp*T=T_y\vp(x-y)\in Y(\Rn)$. $\vp\mapsto T*\vp$ is a continuous linear operator in $Y(\Rn)$.
\end{lemma}

\Proof The first claim follows from Theorem \ref{t1}, (3), the fact that $Y(\Rn)$ is closed under convolution and, finally, from the representation in Theorem \ref{t1}, (2). The second is then easily shown or follows from the Closed Graph Theorem. \qed

This shows part of the following theorem.
\begin{theorem}\label{t2} For an operator $L\in L(Y(\Rn))$ the following are equivalent:
\begin{enumerate}
\item $L$ commutes with translations.
\item There is $T\in\OY$ such that $L\vp = T*\vp$ for all $\vp\in Y(\Rn)$.
\end{enumerate}
\end{theorem}

\Proof (2) $\Rightarrow$ (1) is clear, we have to show the converse. We define $T\in Y(\Rn)'$ by $T\vp:= (L\vp)(0)$. Then by standard arguments  we have $(L\vp)(x)=L(\vp(\cdot+x))(0)=T_y\vp(y+x)=\check{T}\vp(x-y)=(\check{T}*\vp)(x)$. Due to Theorem \ref{t1}, (3) we have $T\in\OY$, hence also $\check{T}\in\OY$ .  \qed

The dual situation is a bit more complicated, since existence of $T*S$ and commutivity is not a priori clear. We define:
$$\Oy:=\{f\in\CN\,:\,\exists k\in\N_0\, \forall \alpha\in\N_0^d: \,\sup_{x\in\Rn} |f^{(\alpha)}(x)|\,e^{-k|x|}<\infty\}.$$
Equipped with its natural locally convex topology $\Oy$ is the inductive limit of a sequence of \F spaces, that is, an (LF)-space and we have

\begin{lemma}\label{l2} $\OY$ is the dual space of $\Oy$. For $S\in Y(\Rn)'$ the map $\vp\mapsto S_y\vp(x+y)$ is a continuous linear map from $Y(\Rn)$ to $\Oy$.
\end{lemma}

\Proof The first part by use of a standard argument using Theorem \ref{t1}, (2). For the second part we estimate
\begin{eqnarray*}
|\partial^\alpha S_y\vp(x+y)|&=&|S_y\vp^{(\alpha)}(x+y)|\le C \sup_{y,|\beta|\le k}|\vp^{(\alpha+\beta)}(x+y)e^{k|y|}\\
&\le& C\,e^{k|y|} \sup_{\xi,|\gamma|\le m} |\vp^{(\gamma)}(\xi)|e^{k|\xi|}
\end{eqnarray*}
with $C$ and $k$ depending on $S$ and $m=k+|\alpha|$. \qed

We obtain an analogue to Lemma \ref{l1}.

\begin{lemma}\label{l3} If $T\in\OY$ and $S\in Y(\Rn)'$ then both $T*S\in Y(\Rn)'$ and $S*T\in Y(\Rn)'$ are defined and equal. $S\mapsto T*S$ is a continuous linear operator in $Y(\Rn)'$.
\end{lemma}

\Proof The existence of $S*T$ follows from Theorem \ref{t1}, (3), the existence of $T*S$ from Lemma \ref{l3}. $(T*S)\vp=(S*T)\vp$ for $\vp\in\Di$ equality follows by direct calculation by use of Theorem \ref{t1}, (2). The continuity of $S\mapsto S*T$ is obvious. \qed

\begin{theorem}\label{t3} For an operator $L\in L(Y(\Rn)')$ the following are equivalent:
\begin{enumerate}
\item $L$ commutes with translations.
\item There is $T\in\OY$ such that $L(S) = T*S$ for all $S\in Y(\Rn)'$.
\end{enumerate}
\end{theorem}

\Proof (2) $\Rightarrow$ (1) is clear, we have to show the converse. The transpose $L^*\in L(Y(\Rn))$ also commutes with translation. Note that $Y(\Rn)$ is Montel, hence reflexive. Because of Theorem \ref{t2}, Proof, there is $T\in\OY$ such that $(L^*(\vp))(x)=T_y\vp(x+y)$. So for $S\in Y(\Rn)'$ we obtain $\langle L(S),\vp\rangle=\langle S,L^*(\vp)\rangle=S_x(T_y \vp(x+y))= ((S*T)\vp)(x)$. \qed

\section{Hadamard operators on $\cS'(Q)$} \label{s3}

Let $L$ be a Hadamard operator on $\cS'(Q)$, that is an operator which admits all monomials as eigen-functions. We need some preparations, cf. Section 1 in \cite{VDprime1}. For $a\in Q$ we define the dilation $D_a\in L(\Sipq)$ by
$$(D_a T)\vp:= T_x\left(\frac{1}{a_1\dots a_d}\,\vp\Big(\frac{x}{a}\Big )\right)$$
for $T\in\Sipq$ and $\vp\in\Siq$. By direct verification we see that $D_a\xi^\alpha=a^\alpha\xi^\alpha$.

Like in \cite[Lemma 1.1]{VDprime1} we obtain that $L$ commutes with dilations, that is, $D_a\circ L=L\circ D_a$ for all $a\in Q$.

We set $M=L^*\in L(\Siq)$ and obtain like in \cite[Lemma 1.3]{VDprime1}  that $M$ commutes with with dilations, that is,
$$M_\xi(\vp(\eta\xi))[x]=(M\vp)(\eta x)$$
for all $\vp\in \Siq$ and $\eta\in Q$.

For $\vp\in\Siq$ we define now
$$T\vp= (M\vp)(\be)=(L\delta_\be)(\vp).$$
Then $T\in\Sipq$ and for all $\eta\in Q$ we have
\begin{equation}\label{e1}(M\vp)(\eta)=T_\xi \vp(\eta\xi).\end{equation}
We have to determine the set of distributions in $T\in\Sipq$ such that
\begin{equation}\label{e2} T_\xi \vp(\cdot\,\xi)\in\Siq \text{ for all } \vp\in \Siq. \end{equation}

For $\widetilde{T}=C_\Log^*(T)$ the condition (\ref{e2}) is equivalent to
\begin{equation}\label{e3} \widetilde{T}_\xi \psi(\cdot+\xi)\in Y(\Rn) \text{ for all } \psi\in Y(\Rn) \end{equation}
which, by Theorem \ref{t1}, is equivalent to $\widetilde{T}\in\OY$.

In analogy to \cite[Definition 3]{VDprime1} we define the space $\OHq$ of \em $\theta$-rapidly decreasing distributions \rm on $Q$.

\begin{definition}\label{d2} $T\in \OHq$ if for any $k$ there are finitely many functions $t_\beta$ such that $(|x|^{2k}+|x|^{-2k})t_\beta\in L_\infty(Q)$ and such that $T=\sum_\beta \theta^\beta t_\beta$.
\end{definition}
 By use of the description in Theorem \ref{t1}, (2), we obtain:

 \begin{lemma} \label{l4} $C_\Exp^*(\OY)=\OHq$.
 \end{lemma}

 Hence we obtain the following translation of Theorem \ref{t1}:
 \begin{theorem}\label{t4}For $T\in\Dip$ the following are equivalent:
\begin{enumerate}
\item $T\in\OHq$
\item For any $k$ there are finitely many functions $t_\beta$ such that $(|x|^{2k}+|x|^{-2k})t_\beta\in L_\infty(Q)$ and such that $T=\sum_\beta \theta^\beta t_\beta$
\item $T\in \cS'(Q)$ and $T_x \vp(xy)\in \cS(Q)$ for all $\vp\in \cS(Q)$.
\item $f(y)=T_x \vp(xy)$ is a rapidly decreasing continuous function (that is $\sup |f(y)|\\ (|y|^{2k}+|y|^{-2k})<\infty$ for all $k$) for all $\vp\in\cD(Q)$.
\item $((|x|^{2k}+|x|^{-2k})T)\star\vp$ is a continuous bounded function for every $k$ and $\vp\in\cD(Q)$.
\end{enumerate}
\end{theorem}

We have obtained the following.
\begin{theorem}\label{t6} $L$ Hadamard operator on $\cS'(Q)$ if and only if there is $T\in\OHq$ such that $L(S)=S\star T$ for all $T\in \cS'(Q)$. \end{theorem}

Here $\langle S\star T,\vp\rangle = S_x(T_y \vp(xy)$ for all $\vp\in\cS(Q)$.

\section{Hadamard operators on $\Sip$} \label{s4}

Let now $L$ be a Hadamard operator on $\Sip$ and $M=L^*\in L(\Si)$ and obtain like in \cite[Lemma 1.3]{VDprime1}  that $M$ commutes with with dilations, that is,
\begin{equation}\label{e6} M_\xi(\vp(\eta\xi))[x]=(M\vp)(\eta x)\end{equation}
for all $\vp\in \Si$ and $\eta\in \NZ$.

For $\vp\in\Si$ we define now
\begin{equation}\label{e9} T\vp= (M\vp)(\be)=(L\delta_\be)(\vp). \end{equation}
Then $T\in\Sip$ and for all $\eta\in \NZ$ we have
\begin{equation}\label{e4}(M\vp)(\eta)=T_\xi \vp(\eta\xi).\end{equation}
We have to determine the set of distributions in $T\in\Sip$ such that
$T_\xi \vp(\cdot\,\xi)$, $\xi\in\NZ$, extends to a function in $\Si$ for all $\vp\in\Si$.

We want to use the results of Section \ref{s3}. We denote by $H_j$, $j=1,..,d$, the coordinate hyperplanes and set $Z_0=\bigcup_j H_j$.
$$\cS(\NZ)=\{\vp\in\Si\,:\, \vp \text{ flat on } Z_0\}.$$
We will show that $M(\cS(\NZ))\subset \cS(\NZ)$. For that it suffices to show that $L(\cS'(Z_0))\subset \cS'(Z_0)$. Here $\cS'(Z_0)$ denotes the temperate distributions with support in $Z_0$.

By $\cF$ we denote the Fourier transform and remark that  for all $j$
\begin{equation}\label{e7}\theta_j\circ \cF=\cF\circ \theta_j^*, \quad \theta_j^*\circ\cF=\cF\circ\theta_j.\end{equation}

We set $\tL=\cF\circ L \circ \cF^{-1}$ and since $\theta_j^*$ commutes with $L$ we conclude by use of (\ref{e7}) that $\tL$ commutes with $\theta_j$ for all $j$. By straightforward calculation we see that $\tL$ commutes with all reflections.  By Lemma \ref{l9} this implies that $\tL$ is a Hadamard operator. We have $L=\cF^{-1}\circ\tL\circ\cF$.

We have $\cF (\delta^{(\alpha)})= i^\alpha (2\pi)^{-d/2} x^\alpha$, hence $\tL(\cF \delta^{(\alpha)})=\widetilde{m}_\alpha \cF \delta^{(\alpha)}$. Finally we obtain
\begin{equation} \label{e8} L(\delta^{(\alpha)})=\widetilde{m}_\alpha \delta^{(\alpha)},\end{equation}
where
$\tL(x^\alpha)=\widetilde{m}_\alpha x^\alpha$.

\sc Example: \rm $L=\theta$ then $\tL = \theta^*=-\theta-1$.  Since $\widetilde{m}_k =-k-1$ we obtain $\theta \delta^{(k)}= (-k-1) \delta^{(k)}$ which, of course, can be verified by direct calculation.

In fact, we will need this result only for $d=1$. We set $x=(x_1,x')$, $x'=(x_2,..,x_d)$ and consider distributions of the form $T_\alpha=\delta^{(\alpha)}(x_1)\otimes S(x')$, $S\in\cS'(\R^{d-1})$.

We fix $\alpha'=(\alpha_2,..,\alpha_d)$ and $\psi\in\cD(\R^{d-1})$. For $T\in\Sip$ and $\vp\in\si$ we set
$(R_\psi T)\vp:= T(\vp(x_1)\psi(x')$.  This defines a map $R_\psi:\Sip\to\sip$.

For $U\in\sip$ we set $L_1(U):=(R_\psi\circ L)(U\otimes x^{\alpha'})\in \sip$. We obtain for $\alpha\in\N_0$ and $\hat{\alpha} = (\alpha,\alpha')$
$$ L_1(x^\alpha) =R_\psi(L x^{\hat{\alpha}})=R_\psi(m_{\hat{\alpha}}x^{\hat{\alpha}})= m_{\hat{\alpha}}\int \xi^{\alpha'}\psi(\xi)\,d\xi\;  x^\alpha.$$
Hence $L_1$ is a Hadamard operator on $\sip$ and, by (\ref{e8}), $\delta^{(\alpha)}$ is an eigenvector of $L_1$.
This means $L_1(\delta^{(\alpha)}) = \mu_\alpha \delta^{(\alpha)}$, hence $(-1)^\alpha \mu_\alpha \vp^{(\alpha)}(0) = T(\vp(x_1),\psi(x'))$
for all $\vp\in\si$ where $T=L(\delta^{(\alpha)}\otimes x^{\alpha'})$.

We choose $\chi\in\di$ , $\chi=0$ in a neighborhood of $0$, and set $\vp_\alpha(x)= \frac{x^\alpha}{\alpha!}\,\chi(x)$.
Then $\mu_\alpha=(-1)^\alpha \, T(\vp_\alpha(x_1),\psi(x'))$. Setting $\mu_\alpha(\psi)=T(\vp_\alpha(x_1),\psi(x'))$ we obtain a distribution $\mu_\alpha\in\cS'(\R^{d-1})$ such that
$$L(\delta^{(\alpha)}\otimes x^{\alpha'}) = \delta^{(\alpha)}\otimes \mu_\alpha.$$

We fix $\alpha\in\N_0$ and we have shown, that $x^{\alpha'}\in \{S\in\cS'(\R^{d-1})\,:\, L(T_\alpha)\in \delta^{(\alpha)}\otimes \cS(\R^{d-1})\}$ for all $\alpha'\in\N_0^{d-1}$. Since this set is a closed linear subspace of $\cS'(\R^{d-1})$ we have shown: $L(T_\alpha)\in \delta^{(\alpha)}\otimes \cS(\R^{d-1}) $ for all $S\in\cS(\R^{d-1})$.

Distributions $T\in \cS'(H_1)$ have the form
$$T=\sum_{\alpha=0}^m \delta^{(\alpha)}(x_1)\otimes S_\alpha(x')$$
(cf. \cite[Chap III, Th\'eor\`eme XXXVI]{LSI}). So we have shown  $L(\cS'(H_1))\subset \cS'(H_1)$. By an analogous argument this holds also for $H_j$, $j=2,..,d$.

Since $\cS'(Z_0)=\sum_{j=1}^d \cS'(H_j)$ (see \cite[Lemma 3.3]{VEs}) we have shown:
\begin{lemma}\label{l7} $L(\cS'(Z_0))\subset \cS'(Z_0)$.
\end{lemma}
As an immediate consequence we obtain:
\begin{proposition}\label{t5} $M(\cS(\NZ))\subset \cS(\NZ)$.
\end{proposition}

We put $M_+ (\vp)= M(\vp)|_{Q}$ for $\vp\in\cS(Q)$. Then $M_+\in L(\cS(Q))$ and $L_+:=M_+^*\in L(\cS'(Q))$ is a Hadamard operator.
From Theorem \ref{t4} we get $T_+=L_+(\be)=T|_{\cS(Q)}\in \OHq$.

Clearly we can do this for all quadrants $Q_e=\{x\,:\,ex\in Q\}$. We set $M_e (\vp)= M(\vp)|_{Q}$ for $\vp\in \cS(Q_e)$ and
$T_e(\vp)=(M_e\vp)(\be)$ for $\vp\in\cS(Q_e)$. By the same arguments as before we obtain that $T_e\in \cO_H'(Q_e)$ (defined in obvious analogy).

In analogy to
Definition \ref{d2} we define the space of distributions on $\Rn$, which are \em $\theta$-rapidly decreasing \rm in infinity and at the coordinate hyperplanes.
\begin{definition}\label{d3} $T\in \OHNZ$ if for any $k$ there are finitely many functions $t_\beta$ such that $(|x|^{2k}+|x|^{-2k})t_\beta\in L_\infty(\Rn)$ and such that $T=\sum_\beta \theta^\beta t_\beta$.
\end{definition}

Then we have for $T$ as defined in (\ref{e9}):

\begin{lemma}\label{t7} $\widetilde{T}:=\sum_e T_e\in\OHNZ$ and $T|_{\cS(\NZ)}=\widetilde{T}$.
\end{lemma}

For $T\in\OHNZ$ and $\vp\in\Si$ we define $(M_T\vp)(x) =T_\xi\vp(\xi x)$ which is defined for all $x\in\Rn$.

\begin{lemma}\label{l8} $M_T$ is a continuous linear operator in $\Si$, $L_T:=M_T^*$ is a Hadamard operator.
\end{lemma}

\Proof We have to estimate $\Psi(x):=x^\gamma (M_T\vp)^{(\alpha)}(x)$. We first recall that $\theta^*_\xi \vp(\xi x) = (\theta^* \vp)(\xi x)$ and $(\theta^*)^\beta \xi^\alpha \vp(\xi) = \xi^\alpha \sum_{\nu\le\beta} p_\nu(\xi)\, \vp^{(\nu)}(\xi)=: \xi^\alpha \psi(\xi)$, where the $p_\nu$ are polynomials.

 We choose $k=k(\alpha-\gamma)$ large enough and obtain
 \begin{eqnarray*} \Psi(x) &=&
 x^\gamma \int \theta^\beta \tau_\beta(\xi)\, \xi^\alpha\,  \vp^{(\alpha)}(\xi x)\,d\xi = x^{\gamma-\alpha}\int \theta^\beta \tau_\beta(\xi) (\xi x)^\alpha  \vp^{(\alpha)}(\xi x)\,d\xi\\
 &=& x^{\gamma-\alpha} \int \tau_{\beta}(\xi)\, (\xi x)^\alpha\, \psi(\xi x)\, d\xi
 = \int \tau_\beta(\xi)\, \xi ^{\alpha-\gamma}\, (\xi x)^\gamma \psi(\xi x) d\xi
 \end{eqnarray*}
 and therefore
 $$\|\Psi\|_{\gamma,\alpha}\le \Big( \int|\tau_\beta(\xi)\,\xi^{\alpha-\gamma}|\, d\xi \Big) \|\vp\|.$$
 Here  $\|\vp\|:= \sup_x |x^\gamma \psi(x)|$ is a continuous semi-norm on $\Si$. This shows the first part of the claim.

 For the second part we have to study $\int x^\gamma (M_T \vp)(x) dx$. We obtain
 \begin{eqnarray*}
\int x^\gamma \Big( \int \theta^\beta \tau_\beta(\xi) \vp(\xi x)\,d\xi\Big)\,dx &=&
\int \theta^\beta \tau_\beta(\xi)\xi^{-\gamma} \Big( \int (\xi x)^\gamma \vp(\xi x)dx \Big)\,d\xi\\ &=&
 \Big( \int \theta^\beta \tau_\beta(\xi) \xi^{-\gamma-\be} d \xi\Big) \int x^\gamma \vp(x) dx\\
 &=&  (-\gamma-\be)^\beta \Big( \int  \tau_\beta(\xi) \xi^{-\gamma-\be} d \xi\Big) \int x^\gamma \vp(x) dx\\
 &=& \int (m_\gamma x^\gamma) \vp(x) dx.
 \end{eqnarray*}
 We have shown that $L_T x^\gamma = m_\gamma x^\gamma$ and this completes the proof. \qed

For $S\in\Sip$, $T\in\OHNZ$  we define $S\star T\in L(\Sip)$ by
$$(S\star T)(\vp)=  S_x(T_\xi\vp(\xi x)) \text{ for all }\vp\in\Si.$$
The following is the main result of this paper.

\begin{theorem}\label{t8}
1.\,For every $T\in\OHNZ$ the map $S\mapsto S\star T$ is a Hadamard operator on $\Sip$.\\
2.\,For every Hadamard operator $L$ on $\Sip$ there is $T\in\OHNZ$ such that $L(S)=S\star T$ for all $S\in\Sip$.
\end{theorem}

\Proof The first part is Lemma \ref{l8}. For the second part we choose $T=L(\delta_\be)$, then $(M\vp)(\eta)=T_\xi \vp(\eta\xi)$ for all $\eta\in\NZ$ (see (\ref{e4})). This equation is true for all $\eta\in\Rn$ if $\vp\in\cS(\NZ)$. By Lemma \ref{t7} there is $\widetilde{T} \in\OH$ such that $T\vp = \widetilde{T}\vp$ for $\vp\in \cS(\NZ)$. This means that $L_{\widetilde{T}}(S)=S\star \widetilde{T}$ defines a Hadamard operator and $L(S)\vp=L_{\widetilde{T}}(\vp)$ for all $\vp\in\cS(\NZ)$. Therefore $L-L_{\widetilde{T}}$ is a Hadamard operator
such that $(L-L_{\widetilde{T}})S$ vanishes on $\cS(\NZ)$ for all $S\in\Sip$, hence $(L-L_{\widetilde{T}})x^\alpha=0$ for all $\alpha$ and therefore $L-L_{\widetilde{T}}=0$. Finally we have $T=L(\delta_\be)=L_{\widetilde{T}}(\delta_\be)=\widetilde{T}$. Therefore we have $L(S)=S\star T$ for all $S\in\Sip$. \qed

\section{Final remarks} \label{s5}

In \cite{VDprime1} the Hadamard operators in $\Dip$ were characterized. We can express the Main Theorem of \cite{VDprime1} in the following way:

\begin{theorem}\label{t9}
The Hadamard operators on $\Dip$ are the operators of the form $S\mapsto S\star T$ where $T\in\OHNZ$ and $\supp T$ has positive distance to the coordinate hyperplanes.
\end{theorem}

This follows from the fact that for a distribution $T\in\Dip$, the support of which has positive distance of the coordinate hyperplanes, the conditions $T\in\OH$ and $T\in\OHNZ$ coincide. For the definition of $\OH$ see \cite[Definition 3]{VDprime1}.

This implies:
\begin{corollary}\label{c1} Every Hadamard operator on $\Dip$ maps $\Sip$ into $\Sip$.
\end{corollary}

By $\sigma(x)=\prod_j \frac{x_j}{|x_j|}$ we denote the signum of $x$. For $\alpha\in\N_0^d$ and $T=\theta^\beta \tau_\beta$ with $(|x|^{2k}+|x|^{-2k})t_\beta\in L_\infty(\Rn)$ and $k$ large enough we define
$$T\Big(\frac{\sigma(x)}{x^{\alpha+\be}}\Big)= \int \tau_\beta(x) (\theta^*)^\beta\frac{\sigma(x)}{x^{\alpha+\be}}dx.$$
Therefore, using a proper representation, we can define $T\big(\frac{\sigma(x)}{x^{\alpha+\be}}\big)$ for any $T\in\OHNZ$. The definition does not depend on the representation, as the following result shows.

\begin{theorem}\label{t10} If $T\in\OHNZ$ and $L(S):= S\star T$ the related Hadamard operator on $\Sip$,  then the eigenvalues of
$L$ with respect to $x^\alpha$ are $m_\alpha=T\big(\frac{\sigma(x)}{x^{\alpha+\be}}\big)$.
\end{theorem}

The \bf Proof \rm is the same as the proof of Theorem 4.2 in $\cite{VDprime1}$. In a remark after the proof there it is pointed out that it holds in a very general context.

We could also, in analogy to Section \ref{s1}, define
$$\Ohnz:=\{f\in C^\infty(\NZ)\,:\,\exists k\in\N_0\, \forall \alpha\in\N_0^d: \,\sup_{x\in\NZ} |\theta^{(\alpha)}f(x)|(|x|^{2k}+|x|^{-2k})<+\infty\}.$$
Equipped with its natural locally convex topology $\Ohnz$ is the inductive limit of a sequence of \F spaces, that is, an (LF)-space and we have
\begin{lemma}\label{l10} $\OHNZ$ is the dual space of $\Ohnz$.
\end{lemma}
This can be derived from Lemma \ref{l2} by use of the exponential transformation applied to all quadrants, or by direct verification.
In this setting the term   $T\big(\frac{\sigma(x)}{x^{\alpha+\be}}\big)$ is properly defined.

\vspace{1cm}

\noindent Bergische Universit\"{a}t Wuppertal,
\newline Dept. of Math., Gau\ss -Str. 20,
\newline D-42119 Wuppertal, Germany
\newline e-mail: dvogt@math.uni-wuppertal.de

\end{document}